\documentclass[amssymb,aps,pra,showpacs,showkeys,reprint,floatfix]{revtex4-1}
\usepackage{graphicx}
\usepackage{amsmath}
\DeclareMathOperator{\lcm}{lcm}

\begin{document}

\title[Circumscribed Regular Polygons]{Tightly Circumscribed Regular Polygons}

\author{Richard J. Mathar}
\homepage{http://www.mpia.de/~mathar}
\email{mathar@mpia.de}
\affiliation{Hoeschstr. 7, 52372 Kreuzau, Germany}

\date{\today}

\begin{abstract}
A regular polygon circumscribing  another regular polygon (with a different side number)
may be tightened to minimize the difference of both areas. The manuscripts computes the
optimum result under the restriction that both polygons are concentric, and obtains
limits if the process is repeated in a two-dimensional Babuschka-doll fashion with side numbers
increasing or decreasing by one or stepping through the prime numbers.
The new aspect compared to the circumscription discussed in the literature so far is that further
squeezing of the outer polygon is possible as we drop the requirement of
drawing intermediate spacing circles between the polygon pairs.
\end{abstract}
\maketitle

\section{Introduction}
\subsection{Notation}

A regular $n$-gon is drawn with $n$ edges of some common side length $s_n$.
The perimeter is $ns_n$.
The incircle with inradius $r_n^{(i)}$ touches each edge at the mid point.
Each edge covers an angle of 
\begin{equation}
\phi_n=\frac{2\pi}{n}
\end{equation}
if viewed from the incircle center.
The area of the polygon is comprised of
$n$ rotated copies of an isosceles
triangle in which
the short edge has length $s_n$, facing the angle
$\phi_n$, and the two other edges have 
length $r_n^{(o)}$, which also is the radius of the circumcircle.
This isosceles triangle might be sliced
into two symmetric rectilinear triangles by drawing a line (apothem)
from its base center to the midpoint of the polygon's incircle;
the definition of the
tangent and sine functions in these yield
\begin{equation}
\tan\frac{\phi_n}{2} = \frac{s_n/2}{r_n^{(i)}},
\end{equation}
\begin{equation}
\sin\frac{\phi_n}{2} = \frac{s_n/2}{r_n^{(o)}},
\end{equation}
and therefore
\begin{equation}
r_n^{(i)} = \frac{s_n}{2\tan\frac{\phi_n}{2}},
\label{eq.rni}
\end{equation}
\begin{equation}
r_n^{(o)} = \frac{s_n}{2\sin\frac{\phi_n}{2}}.
\label{eq.rno}
\end{equation}
As a bridge between two-dimensional geometry and numerical algebra,
we define the \emph{standard} position of the polygon in the
Cartesian $(x,y)$ plane by mapping $x$ and $y$ to the real and imaginary
part of the complex plane, placing the vertices
labeled $j=0,1,\ldots n-1$
counter-clock-wise
at the coordinates
\begin{equation}
x+iy = r_n^{(o)} e^{2\pi j i/n},
\label{eq.std}
\end{equation}
where $i\equiv \sqrt{-1}$ is the imaginary unit.
Edges/sides are also enumerated from $0$ to $n-1$ by calling the
smaller of the two vertex labels that are joined.

\subsection{Tight Circumscription}\label{sec.hist}
The circumscription of a regular $n$-gon by a
regular $m$-gon has been constructed earlier by an elementary step
drawing the $n$-gon, its circumcircle, declaring this circle to
be also the incircle of the $m$-gon with $r_m^{(i)}=r_n^{(o)}$,
and drawing the $m$-gon around its incircle
\cite[p. 428]{Finch}\cite[p. 2300]{WeissteinCRC}\cite[A051762]{EIS}.

The theme of this manuscript is to drop the requirement of
equating the two circles and to search for smaller 
circumscribing regular $m$-gons in the extended
range $r_n^{(i)}\le r_m^{(i)}\le r_n^{(o)}$.

This requires positions where some edges of the $m$-gon cut through
the circumcircle of the $n$-gon. The tighter solution, however, may
exist only within a restricted range of (relative) orientations
of the polygons. The manuscript works out a full representation
of the smallest circumscribing polygons, using the ratio $r_m^{(o)}/r_n^{(o)}$
as a figure of merit.

In overview, the achievable size ratios are calculated in Section \ref{sec.std}
for concentric polygon pairs aligned such that the common center,
a vertex of the inner polygon and a vertex of the outer polygon are collinear
(standard positions).
In Section \ref{sec.rot} further size reductions of the outer polygon are
found if the side number $n$ of the inner polygon is even and the
outer polygon is turned around the common center.
A cursory outlook in Section \ref{sec.trans} shows that shifting the center
of the outer polygon away from the center of the inner polygon may
define even smaller circumscribing $m$-gons.

\section{Concentric Standard Placements} \label{sec.std}
\subsection{General side numbers}

The definition of circumscription implies
that the outer $m$-gon must stay further away
 from the origin
than the inner $n$-gon
at all viewing directions;
at one or more points of contact, both polygons have
the same distance $r_n^{(o)}$ to the origin.

\begin{figure}[hbt]
\includegraphics[scale=0.45]{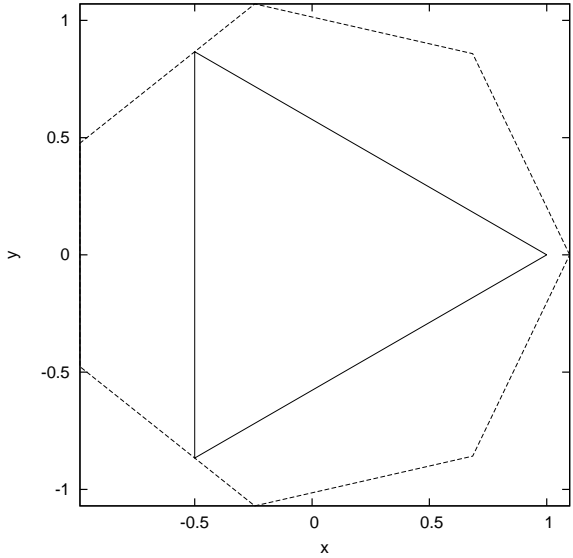}
\caption[]{
Circumscribing a 3-gon by a 7-gon.
}
\label{fig.Rpol37}
\end{figure}

Finding the smallest $m$-gon for small side numbers like
Fig.\ \ref{fig.Rpol37} works as follows.
Any point on the side of the outer polygon of unit radius has a position
$\exp(2\pi i j_o/m)+t[
\exp(2\pi i (j_o+1)/m)
-\exp(2\pi i j_o/m) ]$ with parameter $0\le t\le 1$.
Once the vertex number $j_i$ of the inner polygon
and the side number $j_o$ of the outer polygon which it touches are known,
the point of contact between both polygons in the complex
plane solves
\begin{equation}
r_{m}^{(o)}
\{e^{2\pi i j_o/m}+t[
e^{2\pi i (j_o+1)/m}
-e^{2\pi i j_o/m) }]\}
=
r_{n}^{(o)}e^{2\pi i j_i/n}
.
\label{eq.shrink}
\end{equation}
and after division through $r_m^{(0)}e^{2\pi i j_0/m}$
\begin{equation}
1+t[
e^{2\pi i /m}
-1]
=
\frac{r_{n}^{(o)}}{r_{m}^{(o)}}e^{2\pi i (j_i/n-j_o/m)}
.
\label{eq.jijo}
\end{equation}
Real and imaginary part of this equation establish an inhomogeneous $2\times 2$
linear system of equations for the unknown $t$ and $r_n^{(o)}/r_m^{(o)}$:
\begin{eqnarray}
1+t[
\cos(2\pi /m)
-1]
&=&
\frac{r_{n}^{(o)}}{r_{m}^{(o)}}\cos[2\pi (j_i/n-j_o/m)],
\\
t
\sin(2\pi /m)
&=&
\frac{r_{n}^{(o)}}{r_{m}^{(o)}}\sin[2\pi (j_i/n-j_o/m)]
.
\end{eqnarray}
\begin{eqnarray}
\left(
\begin{array}{cc}
\cos[2\pi (j_i/n-j_o/m)] & 1-\cos(2\pi/m) \\
\sin[2\pi (j_i/n-j_o/m)] & -\sin(2\pi/m)
\end{array}
\right)
&&
\nonumber \\
\cdot
\left(
\begin{array}{c}
r_{n}^{(o)}/r_{m}^{(o)} \\
t
\end{array}
\right)
&=&
\left(
\begin{array}{c}
1\\
0
\end{array}
\right)
.
\end{eqnarray}
The solution is obtained with Cramer's rule.
The inverse ratio is
\cite[4.3.16,4,3,35]{AS}
\begin{eqnarray}
r_{m}^{(o)}/r_{n}^{(o)}
&=&
\frac
{
\sin(j_i\phi_n-j_o\phi_m)
-\sin[j_i\phi_n-(j_o+1)\phi_m] 
}
{\sin \phi_m}
\nonumber \\
&=&
\frac
{
\cos[j_i\phi_n-(j_o+1/2)\phi_m]
}
{\cos( \phi_m/2)}
\nonumber \\
&=&
\frac
{
\cos[\frac{\pi}{nm}\{2j_im-(2j_o+1)n\}]
}
{\cos( \pi/m)}
\label{eq.rrat}
.
\end{eqnarray}
Let $d=m-n$ be the difference in the side numbers; then
\begin{equation}
2j_im-(2j_o+1)n = 2j_id-(2j_o-2j_i+1)n
\label{eq.mism}
\end{equation}
in the argument of the cosine
is a ``mismatch'' value in the angular directions.

The geometric interpretation of this equation:
The minimum radius of the outer polygon is determined by the edge $j_o$
that first hits a nearby vertex $j_i$ while shrinking. The relevant index pair is
the one that maximizes  the cosine in the numerator, so the phase angle is steered towards
zero or $2\pi$, equivalent to $j_i\phi_n\approx (j_0+1/2)\phi_m$. This means the relevant
phase angle in the complex plane and viewing direction is where the vertex $j_i$ points near the middle of edge $j_o$. 

Examples:
\begin{itemize}
\item
In Fig.\ \ref{fig.Rpol37} we have set $n=3$, $m=7$, $r_n^{(o)}=1$
and observe that $j_o=2$, $j_i=1$. Eq.\ (\ref{eq.rrat})
obtains
$r_m^{(o)}/r_n^{(o)}=\cos(\pi/21)/\cos(\pi/7) \approx 1.097519$.
\item
In the inner pair of Fig.\ \ref{fig.Rpol3to16} we have set $n=3$, $m=4$,
$r_n^{(o)}=1$
and observe that $j_o=1$, $j_i=1$. Eq.\ (\ref{eq.rrat})
obtains $r_m^{(o)}/r_n^{(o)}= \sqrt{2}\cos(\pi/12)=(1+\sqrt{3})/2
\approx 1.366025$ \cite{Conwayarxiv98,GirstmairAA81}.
\item
In the inner pair of Fig.\ \ref{fig.Rpol3toPri} we have set $n=3$, $m=5$,
$r_n^{(o)}=1$
and observe that $j_o=1$, $j_i=1$. The equation yields
$r_m^{(o)}/r_n^{(o)}=\cos(\pi/15)/\cos(\pi/5)=\frac{3-\sqrt{5}}{4}+\sqrt{3}\sqrt{\frac{5-\surd 5}{8}}
\approx 1.2090569$.
\end{itemize}

Numerical examples of the size ratios of are gathered
in Table \ref{tab.reg}.
[In solutions with interlaced circles---summarized in
Section \ref{sec.hist}---the ratio $r_m^{(o)}/r_n^{(o)}$ always equals
$1/\cos(\phi_m/2)$ derived with Eqs.\ (\ref{eq.rni}) and (\ref{eq.rno}).
This restricted search space would have put constant values down each
column.]

Where $m$ is a multiple of $n$, the table entries equal one.
In these cases one can re-use $n$ vertices of the inner polygon as vertices
of the outer polygon, and  obtains its remaining $m-n$ vertices
by regular subdivision of the angle,
$
\phi(m) = \frac{\phi(n)}{m/n}.
$
The circumradii are the same, $r_m^{(o)}=r_n^{(o)}$, and their ratio equals one.

\begin{table*}
\begin{tabular}{r|rrrrrrrr}
 & 3 & 4 & 5 & 6 & 7 & 8 & 9 & 10 \\ \hline
 3 & 1.00000000 & 1.36602540 & 1.20905693 & 1.00000000 & 1.09751942 & 1.07313218 & 1.00000000 & 1.04570220 \\ 
 4 & 2.00000000 & 1.00000000 & 1.23606798 & 1.15470054 & 1.10991626 & 1.00000000 & 1.06417777 & 1.05146222 \\ 
 5 & 1.95629520 & 1.39680225 & 1.00000000 & 1.14837497 & 1.10544807 & 1.07905555 & 1.06158549 & 1.00000000 \\ 
 6 & 2.00000000 & 1.36602540 & 1.23606798 & 1.00000000 & 1.10991626 & 1.07313218 & 1.06417777 & 1.04570220 \\ 
 7 & 1.97766165 & 1.40532128 & 1.23109193 & 1.15147176 & 1.00000000 & 1.08068940 & 1.06285492 & 1.05040347 \\ 
 8 & 2.00000000 & 1.41421356 & 1.23606798 & 1.15470054 & 1.10991626 & 1.00000000 & 1.06417777 & 1.05146222 \\ 
 9 & 1.87938524 & 1.40883205 & 1.23305698 & 1.13715804 & 1.10853655 & 1.08136200 & 1.00000000 & 1.05082170 \\ 
 10 & 2.00000000 & 1.39680225 & 1.23606798 & 1.14837497 & 1.10991626 & 1.07905555 & 1.06417777 & 1.00000000 \\ 
\end{tabular}
\caption{Size ratios $r_m^{(o)}/r_n^{(o)}$ as a function of the polygon edge
count (rows $n\ge 3$ and columns $m\ge 3$) for the tight regular concentric positions.
}
\label{tab.reg}
\end{table*}

\subsection{Consecutive side numbers}

\begin{figure}[hbt]
\includegraphics[scale=0.45]{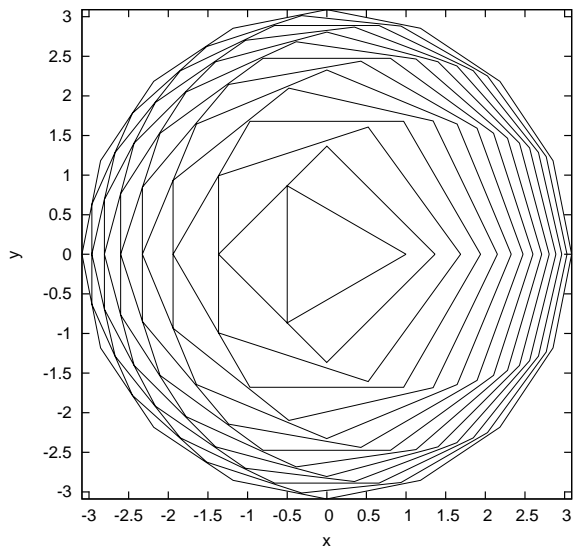}
\caption[]{
Circumscribing a 3-gon by a 4-gon by a 5-gon \ldots up to a 16-gon,
all at  concentric standard positions.
}
\label{fig.Rpol3to16}
\end{figure}

The radial growth of a repeated, possibly infinite nesting
is illustrated in Figure \ref{fig.Rpol3to16}\@.
The circumradius is the partial products of the first upper sub-diagonal of
Table \ref{tab.reg}.

The heuristics is here with $m=n+1$ that
\begin{itemize}
\item
for odd $n$,
the edge $j_o=m/2-1$ hits the vertex $j_i=m/2-1$
with residual mismatch (\ref{eq.mism}) equal to $-1$. [This is the best possible
absolute value because $2j_im$ is even and $(2j_o+1)n$ is odd then.]
So (\ref{eq.rrat}) is
\begin{equation}
r_{n+1}^{(o)}/r_{n}^{(o)}
=
\frac
{
\cos\frac{\pi}{n(n+1)}
}
{\cos \frac{\pi}{n+1}},\quad n\,\mathrm{odd},
\end{equation}
\item
For even $n$ and $m=n+1$, the edge $j_o=n/2$ hits the vertex $j_i=n/2$
on the negative real axis with residual mismatch (\ref{eq.mism}) equal to zero.
So (\ref{eq.rrat}) yields
\begin{equation}
r_{n+1}^{(o)}/r_{n}^{(o)}
=
\frac
{
1
}
{\cos \frac{\pi}{n+1}}, n\,\mathrm{even}.
\end{equation}
\end{itemize}
These two equations constitute the first upper diagonal of Table \ref{tab.reg}\@.
Fencing the polygons up to infinity defines the
limiting
radius as an alternating product of these two factors, 
\begin{eqnarray}
\frac{r_\infty^{(o)}}{r_3^{(o)}}
&=&
\frac{r_4^{(o)}}{r_3^{(o)}}
\times
\frac{r_5^{(o)}}{r_4^{(o)}}
\times
\frac{r_6^{(o)}}{r_5^{(o)}}
\times
\cdots
\nonumber\\
&=&
\prod_{n=3,5,7,\ldots} \frac{ \cos\frac{\pi}{n(n+1)} }{\cos\frac{\pi}{n+1}}
\prod_{n=4,6,8,\ldots} \frac{1}{\cos\frac{\pi}{n+1}}
\nonumber \\
&=&
\frac{\prod_{n=3,5,7,\ldots} \cos\frac{\pi}{n(n+1)}}
{\prod_{n=3}^\infty \cos\frac{\pi}{n+1}}
\nonumber \\
&=&
\frac{1}{2K'}
\prod_{n=3,5,7,\ldots} \cos\frac{\pi}{n(n+1)}
\nonumber \\
&\approx& 4.16674437148793
\label{eq.Couse}
\end{eqnarray}
where we have inserted (\ref{eq.Co})
and
\cite{Finch}\cite[A085365]{EIS}
\begin{equation}
K'\equiv \prod_{n=3}^\infty \cos\frac{\pi}{n}
\approx
0.1149420448532962007.
\label{eq.Kprime}
\end{equation}

\begin{figure}[hbt]
\includegraphics[scale=0.45]{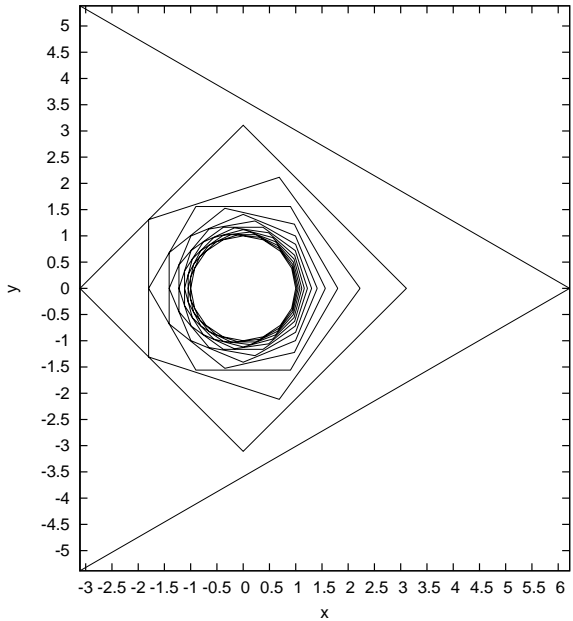}
\caption[]{
Circumscribing a 16-gon by a 15-gon by a 14-gon \ldots down to a 3-gon,
all with the same center.
}
\label{fig.Rpol16to3}
\end{figure}

In Figure \ref{fig.Rpol16to3} the edge count of the circumscribed polygon
is \emph{decreased} from 16 to 3\@. The ratio of the circumradii of the
3-gon and  the 16-gon is $\approx 6.2$ in the image.
The observation is that here
for $m=n-1$ and
\begin{itemize}
\item
$n$ even,
the outer edge $j_o=n/2-1$ and $j_i=n/2$.
The value of (\ref{eq.rrat}) is
\begin{equation}
r_{n-1}^{(o)}/r_{n}^{(o)}
=
\frac
{
1
}
{\cos\frac{\pi}{n-1}}
,\quad n\, \mathrm{even}.
\end{equation}
\item
whereas
for $n$ odd $j_o=(n-3)/2$ and $j_i=(n-1)/2$.
The value of (\ref{eq.rrat}) is
\begin{equation}
r_{n-1}^{(o)}/r_{n}^{(o)}
=
\frac
{
\cos\frac{\pi}{n(n-1)} 
}
{\cos\frac{\pi}{n-1}}
,\quad n\,\mathrm{odd}
\end{equation}
\end{itemize}
The alternating infinite product of these
terms, the finite radius of the free inner region in Figure \ref{fig.Rpol16to3}
if inscribing indefinitely:
\begin{eqnarray}
&& r_3^{(o)}/r_{\infty}^{(o)}
=
\frac{r_3^{(o)}}{r_4^{(o)}}
\times
\frac{r_4^{(o)}}{r_5^{(o)}}
\times
\frac{r_5^{(o)}}{r_6^{(o)}}
\times\cdots
\nonumber \\
&&=
\prod_{n=4,6,8,\ldots}
\frac{r_{n-1}^{(o)}}{r_n^{(o)}}
\times
\frac{r_{n}^{(o)}}{r_{n+1}^{(o)}}
\nonumber \\
&&=
\prod_{n=4,6,8\ldots}
\frac{1}{\cos\frac{\pi}{n-1}}
\times
\frac
{
\cos\frac{\pi}{(n+1)n} 
}
{\cos\frac{\pi}{n}}
\nonumber \\
&& \approx 8.5526818319553
\label{eq.Kstd}
.
\end{eqnarray}
This is slightly smaller than the equivalent
polygon circumscribing constant
$1/K' \approx 8.7000366\ldots$ \cite[A051762]{EIS}\cite[p. 2300]{WeissteinCRC}\cite[p. 428]{Finch,GrimstoneMG64}
by the factor
\begin{equation}
\prod_{n=4,6,8,\ldots} \cos\frac{\pi}{n(n+1)}\approx 0.98306273874458351\ldots
,
\label{eq.cosCorr}
\end{equation}
based on (\ref{eq.Kprime}).
The areas have been smaller
relative to the published
construction with interspersed circles.
The logarithm of the new constant is evaluated in Appendix \ref{sec.logcos}.

\subsection{Prime side numbers}

\begin{figure}[hbt]
\includegraphics[scale=0.45]{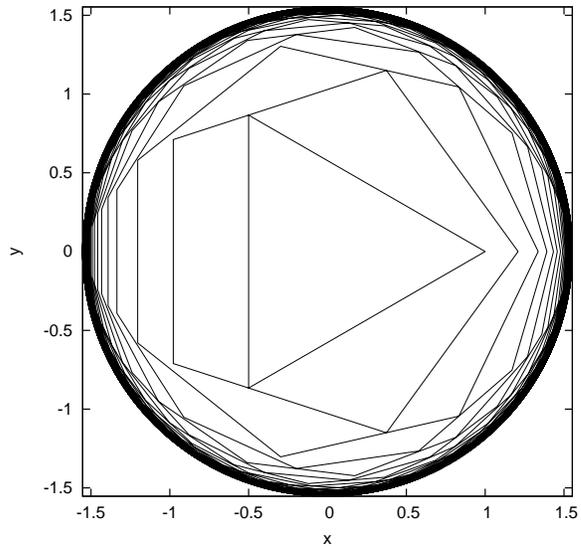}
\caption[]{
Circumscribing a 3-gon by a 5-gon by a 7-gon by a 11-gon etc up to a 541-gon,
all with the same center,
using all odd primes as edge numbers.
}
\label{fig.Rpol3toPri}
\end{figure}

If $m$ is the next prime after $n$, the indices of the vertex of $n$
and edge of $m$ that describes the contact is irregular,
see Fig.\ \ref{fig.Rpol3toPri} and Table \ref{tab.rrpr}.
It is given by the pair $(j_o,j_i)$ which maximizes the value of 
(\ref{eq.rrat}).
The mismatch of (\ref{eq.mism}), $2j_im-(2j_o+1)n$ cannot be nulled for odd primes $m$ and $n$
because $2j_im$ is even and $(2j_o+1)n$ is odd. But the value can
apparently be forced to
$\pm 1$ (where the sign is not important
because this is an argument to the even function of the cosine), as demonstrated in
Table \ref{tab.rrpr}\@.
[The value of Eq.\ (\ref{eq.mism}) is either $+1$ or $-1$
depending on whether the odd number $(2j_o-2j_i+1)n$ is to be
incremented or decremented to reach a multiple of 4, since the
prime gaps $d$ are even and the values of $2j_id$ are multiples of 4.]

Assuming this $\pm 1$ heuristics if always correct, the circumcircle
radius in Figure \ref{fig.Rpol3toPri} grows to
\begin{eqnarray}
\frac{r_\infty^{(o)}}{r_3^{(o)}}
&=&
\frac{\cos\frac{\pi}{nm}}{\cos\frac{\pi}{m}}
=
\frac{\cos\frac{\pi}{3\cdot 5}}{\cos\frac{\pi}{5}}
\times
\frac{\cos\frac{\pi}{5\cdot 7}}{\cos\frac{\pi}{7}}
\times
\frac{\cos\frac{\pi}{7\cdot 11}}{\cos\frac{\pi}{11}}
\times
\cdots
\nonumber \\
&=&
\frac{1}{2}
\frac{1}{\cos\frac{\pi}{3}}
\frac{\cos\frac{\pi}{3\cdot 5}}{\cos\frac{\pi}{5}}
\times
\frac{\cos\frac{\pi}{5\cdot 7}}{\cos\frac{\pi}{7}}
\times
\frac{\cos\frac{\pi}{7\cdot 11}}{\cos\frac{\pi}{11}}
\times
\cdots
\nonumber \\
&=&
\frac{1}{2K_p'}
\prod_{p_j\ge 3} \cos\frac{\pi}{p_jp_{j+1}}
\approx 1.5550895739...
\label{eq.ratPri}
\end{eqnarray}
where 
\begin{equation}
K_p' = \prod_{p=3,5,7,11\ldots} \cos\frac{\pi}{p} \approx 0.3128329
\end{equation}
is Kitson's product over odd primes $p$
\cite{KitsonArxiv06,KitsonMG92}\cite[A131671]{EIS}.
The infinite product of cosines in (\ref{eq.ratPri})
is evaluated in (\ref{eq.Cprim}) and
smaller than unity,
so the equation says that
our construction squeezes the circumradius of the casting prime-sided
regular polygons by more than a factor two compared to Kitson's
variant of construction.

Fig.\ \ref{fig.Rpol3toPri} illustrates why:
A non-zero mismatch angle reflects that
no vertex of the inner polygon
touches a midpoint of a side of the outer polygon;
in consequence the circumradius of the 
inner polygon is larger than the inradius
of the outer polygon for each individual pair of polygons.

\begin{table}
\begin{tabular}{rrrrr}
$n$ & $m$ & $j_i$ & $j_o$ & $2j_im-(2j_o+1)n$ \\
\hline
3 & 5 & 1 & 1 & 1\\
5 & 7 & 1 & 1 & -1\\
7 & 11 & 1 & 1 & 1\\
11 & 13 & 3 & 3 & 1\\
13 & 17 & 5 & 6 & 1\\
17 & 19 & 4 & 4 & -1\\
19 & 23 & 7 & 8 & -1\\
23 & 29 & 2 & 2 & 1\\
29 & 31 & 7 & 7 & -1\\
31 & 37 & 13 & 15 & 1\\
37 & 41 & 14 & 15 & 1\\
41 & 43 & 10 & 10 & -1\\
43 & 47 & 16 & 17 & -1\\
47 & 53 & 4 & 4 & 1\\
53 & 59 & 22 & 24 & -1\\
59 & 61 & 15 & 15 & 1\\
61 & 67 & 5 & 5 & -1\\
67 & 71 & 25 & 26 & -1\\
71 & 73 & 18 & 18 & 1\\
\end{tabular}
\caption[]{
Vertex and edge indices $j_i$ and $j_o$ that maximize (\ref{eq.rrat})
for adjacent primes $n$ and $m$, describing the concentric polygons 
in Fig.\ \ref{fig.Rpol3toPri}.
}
\label{tab.rrpr}
\end{table}

If polygons with sides of odd prime numbers are stacked
in \emph{reverse} order, inscribing a 5-gon in a 3-gon, a 7-gon in the 5-gon,
a 11-gon in the 7-gon etc., there is no substantial modification to
the calculation, because interchanging the values of $n$ and $m$ in
Table \ref{tab.rrpr} appears to  lead again to a list of $\pm 1$ in the mismatches.
Now the ratio of the circumradius of the triangle divided by the
radius of the circular inner hole is
\begin{eqnarray}
\frac{r_3^{(o)}}{r_{\infty}^{(o)}}
=
\cdots \times
\frac{r_7^{(o)}}{r_{11}^{(o)}}
\times
\frac{r_5^{(o)}}{r_7^{(o)}}
\times
\frac{r_3^{(o)}}{r_5^{(o)}}
\nonumber \\
=
\frac{\prod_{p_j\ge 3}\cos\frac{\pi}{p_jp_{j+1}}}
{\prod_{p_j\ge 3} \cos\frac{\pi}{p_j}}
=
\frac{1}{K_p'}
\prod_{p_j\ge 3}\cos\frac{\pi}{p_jp_{j+1}},
\end{eqnarray}
so there is a straight factor of 2 relative to the value
in (\ref{eq.ratPri}).

\section{Concentric, Rotations allowed}\label{sec.rot}

If the outer polygon is rotated by an angle $\alpha_m$ relative
to the standard position (\ref{eq.std}), the vertices move to
\begin{equation}
x+iy = r_m^{(o)} e^{2\pi j i/m+\alpha_m},\quad 0\le j< m.
\label{eq.rot}
\end{equation}

In consequence, all three factors on the left hand side
of (\ref{eq.shrink}) are multiplied by $e^{i\alpha_m}$,
Eq.\ (\ref{eq.jijo}) obtains an additional factor $e^{-i\alpha_m}$ on the right hand side,
and the phase shift
finally enters Eq.\ (\ref{eq.rrat}):
\begin{equation}
\frac{r_m^{(o)}}{r_n^{(o)}}
 = \frac{\cos[j_i\phi_n-(j_o+1/2)\phi_m-\alpha_m]}{\cos(\phi_m/2)}
.
\label{eq.ratWalpha}
\end{equation}

\begin{figure}[hbt]
\includegraphics[scale=0.45]{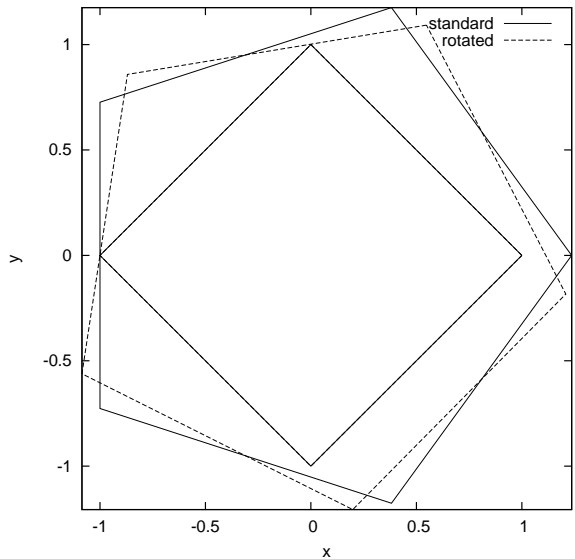}
\caption[]{Tight 5-gon around the 4-gon in the standard placement,
where the angle $\alpha_5$ in (\ref{eq.rot}) is kept at zero, and
another---tighter---solution where $\alpha_5$ is set to $-\pi/20=-9^\circ$
to yield a smaller 5-gon.
The standard solution generates $r_5^{(o)}/r_4^{(o)}\approx 1.236$---see Table \ref{tab.reg}---
whereas the solution allowing rotation yields
$r_5^{(o)}/r_4^{(o)}\approx 1.222$---see Table \ref{tab.rot}.
}
\label{fig.Rpol45a}
\end{figure}

\begin{figure}[hbt]
\includegraphics[scale=0.45]{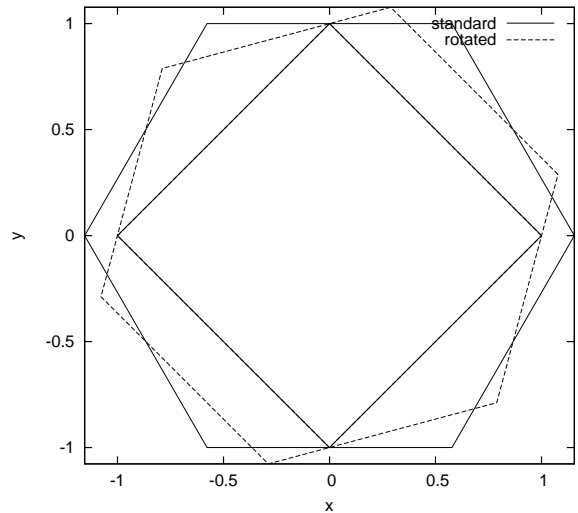}
\caption[]{Tight 6-gon around the 4-gon in the standard placement [where
the angle $\alpha_6$ in (\ref{eq.rot}) is kept at zero], and another tighter
solution where $\alpha_6$ is set to $\pi/12=15^\circ$ to yield a smaller 6-gon.
The standard placement achieves $r_6^{(o)}/r_4^{(o)}\approx 1.154$ according
to Table \ref{tab.reg}, and the version allowing rotation achieves
$r_6^{(o)}/r_4^{(o)}\approx 1.115$ reported in Table \ref{tab.rot}.
}
\label{fig.Rpol46a}
\end{figure}

The tightest solution for fixed $\alpha_m$ is represented by
the pair $(j_i,j_0)$ which maximizes the value of
$r_m^{(o)}/r_n^{(o)}$ and maximizes the value
of the cosine in the numerator
(because $m$ and the denominator are fixed).
Shifts of $\alpha_m$ induce
reduction of some peaks and rises of others in the bi-periodic domain
spanned by the
$j_i$ and $j_o$.
The best solution is obtained
where the value of the cosine becomes degenerate with highest
multiplicity on the grid of the ($j_i$, $j_o$).
In geometrical terms, rotated solutions
seek to maximize the number of contact points between the two polygons,
illustrated in Fig.\ \ref{fig.Rpol45a} and \ref{fig.Rpol46a}.

The interesting  range
is $0\le \alpha_m \le \min(\phi_m,\phi_n)$,
because rotation of the inner polygon by integer multiples of $\phi_n$ or rotation
of the outer polygon by integer multiples of $\phi_m$
leaves the graph invariant.
[Or, formally speaking, changes of $\alpha_m$ modulo $\phi_n$ or modulo $\phi_m$
can be absorbed into resetting the integers $j_i$ or $j_o$ in the numerator.]
The ratio $r_m^{(o)}/r_n^{(o)}$ exercises $m$ and also $n$ periods if $\alpha$ is
turned through a full angle of $2\pi$, and contains therefore $\lcm(m,n)$ periods.
($\lcm$ is the least common multiple of both.)

Because $j_i\phi_n$ and $(j_o+1/2)\phi_m$ have an integer representation if measured in units of
$\pi/(mn)$, because the cosine is a smooth function of its argument,
because the periodicity with respect to $\alpha_m$ means its extremal
values
can only occur at multiples of half the period,
and because $nm=\lcm(n,m)\gcd(n,m)$,
we may encode all relevant angles as
$\alpha_m \equiv s_{n,m}\pi/(mn)$ with integer-valued $s_{n,m}$.
The phase angle in the numerator of (\ref{eq.ratWalpha}) becomes
\begin{equation}
\frac{\pi}{nm}[2j_im-(2j_0+1)n-s_{n,m}]
.
\label{eq.phass}
\end{equation}

Investigation all possible pairs of polygons up to the $88$-gon
leads to the following heuristics:
\begin{itemize}
\item
If $n$ is odd, then $s_{n,m}=0$.
[Interpretation: the mismatch (\ref{eq.mism}) is odd; no vertex
points exactly to the center of an edge. This establishes the following
stability/frustration argument: By the up-down symmetry of the graph,
infinitesimal rotation of the inner polygon requires pushing at least that edge
of the outer polygon outwards, which necessarily growth in size instead
of shrinking as requested.]
This implies that neither the polygon pair in Figure \ref{fig.Rpol37}
nor the cascaded stack with the primal edges numbers
in Figure \ref{fig.Rpol3toPri} can be compressed by adding rotations.
\item
Periodicity:
$s_{n,m}=s_{n,m+n}$.
This seems to be a consequence
of the modular property mentioned above; a change of $m$ by a
multiple of $n$ is absorbed by modifying $j_i$ or $j_o$ by
integer units. The $\gcd(n,m)$ (the period length of the cosine)
is also preserved. Both aspects combined seem to freeze
the number of the contacts between the two polygons.
\item
$s_{n,n}=0$.
If the edge numbers are equal, the circumscribed polygon is a copy
of the inscribed polygon.
\item
If $n$ is even,
\begin{itemize}
\item
$s_{n,n/2}=n/2$. This says that an outer polygon with
half as many vertices as the inner polygon may be constructed
by outwards extension of one over the other edge
of the inner polygon
(which requires a rotation by half of the
angle $\phi_m$ relative to the standard position).
This achieves $r_m^{(i)}=r_n^{(i)}$.
\item
A half period exists with palindromic symmetry:
$s(n,n/2+k)=s(n,n/2-k)$.
Reason: The mirror symmetry of the standard placement leads to
equivalent solutions if the outer polygon
is rotated either clockwise or counter-clockwise.
Solutions are \emph{even} functions of $\alpha_m$, so sign
flips of $s_{n,m}$ are irrelevant. The half period then results
from a general property of (Fourier series of)
periodic even functions.
\end{itemize}
\end{itemize}

Consuming these rules, we need to tabulate the $s_{n,m}$
only in the triangle of even $n$ with $0\le m\le n/2$ for a full coverage.
Then
\begin{itemize}
\item
If $m\le n/2$ is odd and $n$ is even, then
$s(n,m)=\gcd(n/2,m)$.
\item
If $m< n/2$ is even, and
\begin{itemize}
\item
$n$ is two times an odd number, $s_{n,m}=0$.
\item
$n$ is two times an even number, $s_{n,m}=2\gcd(n/2,m/2)$.
This selection is apparently aligning edge 0 of the circumscribing polygon
parallel to edge 0 of
the inscribed polygon with the aim to increase the
number of contacts to a multiple of four, similar to what is observed in Fig.\ \ref{fig.Rpol46a}.
\end{itemize}
\end{itemize}

\begin{figure}[hbt]
\includegraphics[scale=0.46]{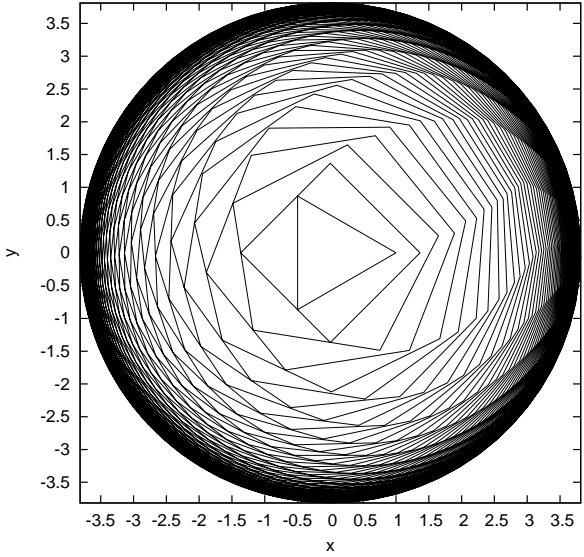}
\caption[]{
Concentric encircling the 3-gon by a 4-gon by a 5-gon and so on
as in
Figure \ref{fig.Rpol3to16}, but minimizing the areas from the 4-gon upwards by
rotating these polygons by variation of $\alpha_m$.
}
\label{fig.Rpol3to16a}
\end{figure}

As an application, the cumulative wrench angle of the vertex direction of the
outermost regular polygon in Figure \ref{fig.Rpol3to16a} relative to its
position in Figure \ref{fig.Rpol3to16}
is calculated as
\begin{eqnarray}
\sum_{n\ge 3,m=n+1} \pi\frac{s_{n,m}}{nm}
=
\pi \sum_{n=2,3,4,\ldots} \frac{s_{2n,2n+1}}{2n(2n+1)}
\nonumber \\
=
\pi \sum_{n=2,3,4,\ldots} \frac{s_{2n,1}}{2n(2n+1)}
=
\pi \sum_{n=2,3,4,\ldots} \frac{1}{2n(2n+1)}
\nonumber \\
=
\pi[\frac{5}{6}-\log(2)]\approx 25.23^\circ.
\end{eqnarray}

Another heuristic observation with this rule is that the absolute value
of the mismatch (\ref{eq.mism}) for $m=n+1$ is kept at 1 if $n$ is odd and at 2 if $n$ is even.
The growth of the radius in Figure \ref{fig.Rpol3to16a} is limited to
the infinite product of terms of the form (\ref{eq.ratWalpha}),
\begin{eqnarray}
\frac{r_{\infty}^{(o)}}
{r_3^{(o)}}
=
\frac{\cos\frac{2\pi}{3\cdot 4}}{\cos\frac{\pi}{4}}
\times \frac{\cos\frac{\pi}{4\cdot 5}}{\cos\frac{\pi}{5}}
\times \frac{\cos\frac{2\pi}{5\cdot 6}}{\cos\frac{\pi}{6}}
\times \frac{\cos\frac{\pi}{6\cdot 7}}{\cos\frac{\pi}{7}}
\times\cdots
\nonumber \\
=
\frac{1}{2K'}
\prod_{k=3,5,7,\ldots}\cos\frac{2\pi}{k(k+1)}
\prod_{k=4,6,8,\ldots}\cos\frac{\pi}{k(k+1)}
\nonumber \\
=
\frac{1}{2K'\cos\frac{\pi}{2\cdot 3}}
C_e
\prod_{k=3,5,7,\ldots}\cos\frac{2\pi}{k(k+1)}
\nonumber \\
=
\frac{1}{\sqrt{3}K'}
C_e
\prod_{k=3,5,7}\cos\frac{2\pi}{k(k+1)}
\approx
3.5809046865583,
\label{eq.roteoInf}
\end{eqnarray}
where
$K'$, $C_e$ and the infinite product
are taken from (\ref{eq.Kprime}), (\ref{eq.Ce}) and (\ref{eq.cos2pi}).
This circumradius including rotations
is considerably smaller than the circumradius
(\ref{eq.Couse}) in the standard positions.

\begin{table*}
\begin{tabular}{r|rrrrrrrr}
 & 3 & 4 & 5 & 6 & 7 & 8 & 9 & 10 \\ \hline
 3 & 1.00000000 & 1.36602540 & 1.20905693 & 1.00000000 & 1.09751942 & 1.07313218 & 1.00000000 & 1.04570220 \\ 
 4 & 1.93716632 & 1.00000000 & 1.22204076 & 1.11535507 & 1.10348396 & 1.00000000 & 1.06044555 & 1.03851698 \\ 
 5 & 1.95629520 & 1.39680225 & 1.00000000 & 1.14837497 & 1.10544807 & 1.07905555 & 1.06158549 & 1.00000000 \\ 
 6 & 1.73205081 & 1.36602540 & 1.22929667 & 1.00000000 & 1.10681271 & 1.07313218 & 1.04801052 & 1.04570220 \\ 
 7 & 1.97766165 & 1.40532128 & 1.23109193 & 1.15147176 & 1.00000000 & 1.08068940 & 1.06285492 & 1.05040347 \\ 
 8 & 1.98422940 & 1.30656296 & 1.23255619 & 1.14559538 & 1.10830702 & 1.00000000 & 1.06324431 & 1.04847492 \\ 
 9 & 1.87938524 & 1.40883205 & 1.23305698 & 1.13715804 & 1.10853655 & 1.08136200 & 1.00000000 & 1.05082170 \\ 
 10 & 1.98904379 & 1.39680225 & 1.17557050 & 1.14837497 & 1.10879865 & 1.07905555 & 1.06352950 & 1.00000000 \\ 
\end{tabular}
\caption{Size ratios $r_m^{(o)}/r_n^{(o)}$ as a function of the polygon edge
counts for rows $n\ge 3$ and columns $m\ge 3$ for the tight
concentric placements, allowing for rotations.
By construction, the elements are not larger than the equivalent
entries in table \ref{tab.reg}.
}
\label{tab.rot}
\end{table*}

\section{Translated Centers}\label{sec.trans}
A glance at
Figure \ref{fig.Rpol3to16} or \ref{fig.Rpol45a} for example
reveals that further compression of the outer polygon
would be possible if either one is shifted sideways,
giving up
the requirement that the two polygons be concentric.

In the complex plane this adds a displacement $z_m^{(o)}$
of the outer polygon as a new parameter to Eq.\ (\ref{eq.shrink}):
\begin{eqnarray}
&&z_m^{(o)}
+
r_{m}^{(o)}
\{e^{2\pi i j_o/m}+t[
e^{2\pi i (j_o+1)/m}
-e^{2\pi i j_o/m) }]\}
\nonumber \\
&&=
r_{n}^{(o)}e^{2\pi i j_i/n}
.
\end{eqnarray}
Assuming that $z_m^{(o)}$ is real-valued (that center shifts are sideways),
this can be written as
\begin{equation}
\frac{z_m^{(o)}}{r_m^{(o)}}\cos(j_o\phi_m)
+
1+t[
\cos(\phi_m)
-1]
=
\frac{r_{n}^{(o)}}{r_m^{(o)}}\cos(j_i\phi_n-j_o\phi_m)
.
\end{equation}
\begin{equation}
-\frac{z_m^{(o)}}{r_m^{(o)}}\sin(j_o\phi_m)
+ t \sin(\phi_m)
=
\frac{r_{n}^{(o)}}{r_m^{(o)}}\sin(j_i\phi_n-j_o\phi_m)
.
\end{equation}
We do not discuss this parameter
space systematically or solutions obtained
by
combined translations and rotations, 
but merely illustrate this aspect by the simplest
examples:
\begin{itemize}
\item
The triangle $n=3$ could touch the quadrangle
$m=4$ in Figure \ref{fig.Rpol3to16}
at its right vertex on the horizontal axis,
at $j_o=j_i=t=0$, which gives
\begin{equation}
\frac{z_m^{(o)}}{r_m^{(o)}}
+ 1
=
\frac{r_{n}^{(o)}}{r_m^{(o)}}
\label{eq.trStd}
\end{equation}
and $0=0$ for the imaginary part.
(This equation and overlap of the two vertices is possible whenever $m>n$.)
The other two vertices of the triangle
would stay glued to the quadrangle's sides,
one point at $j_i=1=j_o$, which is
\begin{equation}
1-t
=
\frac{r_{n}^{(o)}}{r_m^{(o)}}\cos(2\pi/12)
.
\end{equation}
\begin{equation}
-\frac{z_m^{(o)}}{r_m^{(o)}}
+ t
=
\frac{r_{n}^{(o)}}{r_m^{(o)}}\sin(2\pi/12)
.
\end{equation}
The solution to these three linear equations
for the three unknown $t$, $z_m^{(o)}/r_m^{(o)}$ and $r_n^{(o)}/r_m^{(o)}$
is
\begin{eqnarray}
\frac{z_m^{(o)}}{r_m^{(o)}} = 1-\frac{2}{\surd 3}\approx -0.1547005,
\\
\frac{r_n^{(o)}}{r_m^{(o)}}=2(1-\frac{1}{\sqrt{3}})\approx 0.84529946.
\end{eqnarray}
Therefore $r_m^{(o)}/r_n^{(o)}\approx 1.1830127$
which is indeed smaller than $r_4^{(o)}/r_3^{(o)}$ 
in Tables \ref{tab.reg} and \ref{tab.rot}.
\item
For $n=4$, $m=3$, the 3-gon circumscribing the 4-gon in Figure
\ref{fig.Rpol16to3},
the shift leads for the contact on the negative real line where $t=1/2$,
$j_o=(m-1)/2$, $j_i=n/2$ to
\begin{equation}
-\frac{1}{2}\frac{z_m^{(o)}}{r_m^{(o)}}
+\frac{1}{4}=
\frac{1}{2}\frac{r_n^{(o)}}{r_m^{(o)}}
.
\end{equation}

[The generic equation if $n$ is even, $m$ is odd and $n>m$
is
\begin{equation}
-\frac{z_m^{(o)}}{r_m^{(o)}}\sin(\frac{\phi_m}{2})
+\frac{1}{2}\sin(\phi_m)=
\frac{r_n^{(o)}}{r_m^{(o)}}\sin(\frac{\phi_m}{2})
.
\end{equation}
]
Two more equations are established by $j_i=1$, $j_o=0$
if the upper and lower vertex of the quadrangle meets the other
two edges of the triangle,
namely
\begin{eqnarray}
\frac{z_m^{(o)}}{r_m^{(o)}} +1-\frac{3}{2}t=0;
\\
\frac{\sqrt{3}}{2}t=
\frac{r_n^{(o)}}{r_m^{(o)}}
\end{eqnarray}
such that 
\begin{equation}
\frac{r_n^{(o)}}{r_m^{(o)}}=\frac{3}{4}(\sqrt{3}-1)
\approx 1/1.821367
\end{equation}
with center displacement
\begin{equation}
\frac{z_m^{(o)}}{r_m^{(o)}}
=\frac{1}{4}(5-3^{3/2}).
\end{equation}
\item
For the 5-gon circumscribing the 3-gon in Figure
\ref{fig.Rpol3toPri}
the shift leads for the contact on the positive real
axis again to (\ref{eq.trStd})
plus two equations established by $j_i=j_o=1$,
such that 
\begin{eqnarray}
\frac{r_n^{(o)}}{r_m^{(o)}}=
-\frac{2}{3}-\frac{4}{3}\cos\frac{2\pi}{5}
+\frac{2}{3}\cos\frac{2\pi}{15}
+2\cos\frac{4\pi}{15}
\nonumber \\
\approx 1/1.1512750
\end{eqnarray}
with center displacement
\begin{equation}
\frac{z_m^{(o)}}{r_m^{(o)}}
=
-\frac{5}{3}-\frac{4}{3}\cos\frac{2\pi}{5}
+\frac{2}{3}\cos\frac{2\pi}{15}
+2\cos\frac{4\pi}{15}
.
\end{equation}
\end{itemize}

\section{Summary}
We have defined and computed the smallest ratio of the circumradii
of a pair of non-overlapping concentric regular polygons, and have pointed at 
infinite products of cosines that arise if some infinite sets of regular polygons
are nested defined by simple strides in the sets of side numbers.

\appendix

\section{Quenching Factor of the Kepler-Bouwkamp Constant}\label{sec.logcos}
\subsection{Even lower term in the product}

The constant (\ref{eq.cosCorr})  is approached by calculating its logarithm
(and including one more term to put the result into a more general perspective)
\cite{StephensMG79},
\begin{eqnarray}
\log \prod_{n=2,4,6,8,\ldots} \cos\frac{\pi}{n(n+1)}
\nonumber\\
=
\sum_{n=2,4,6,8,\ldots} \log\cos\frac{\pi}{n(n+1)}
\nonumber\\
\approx -0.160923373349205036366901529
\label{eq.numlog}
\end{eqnarray}
via the associated Taylor series \cite[A046991]{EIS}\cite[1.518]{GR}
\begin{equation}
\log\cos\epsilon = -\frac{\epsilon^2}{2}-\frac{\epsilon^4}{12}-\frac{\epsilon^6}{45}-\frac{17\epsilon^8}{2520}
-\frac{31\epsilon^{10}}{14175}
-\frac{691\epsilon^{12}}{935550}
-\cdots
\label{eq.taylogcos}
\end{equation}
as follows:
\begin{eqnarray}
&&-\log \prod_{n=2,4,6,8,\ldots} \cos\frac{\pi}{n(n+1)}
\nonumber\\
=
\sum_{k=1}^\infty
&&
 [
\frac{\pi^2}{2(2k)^2(2k+1)^2}
+\frac{\pi^4}{12(2k)^4(2k+1)^4}
\nonumber\\
&&
+\frac{\pi^6}{45(2k)^6(2k+1)^6}
+\frac{17\pi^8}{2520(2k)^8(2k+1)^8}+\cdots
]
\label{eq.logcost}
.
\end{eqnarray}
Partial fraction decompositions of the individual terms
have the following format
\cite[2.102]{GR}
\cite{MahoneyJCAM9,VellemanAMM109,Xinarxiv04,EusticeAMM86}
\begin{equation}
\frac{1}{n^2(n+1)^2}
=\frac{1}{n^2}+\frac{1}{(n+1)^2}-\frac{2}{n(n+1)},
\end{equation}
\begin{eqnarray}
\frac{1}{n^4(n+1)^4}
=
\frac{1}{n^4}+\frac{1}{(n+1)^4}
-\frac{4}{n^3}+\frac{4}{(n+1)^3}
\nonumber\\
+\frac{10}{n^2}+\frac{10}{(n+1)^2}
-\frac{20}{n(n+1)},
\end{eqnarray}
\begin{equation}
\frac{1}{n^{2s}(n+1)^{2s}}
=
\sum_{t=1}^{2s}\
\binom{4s-t-1}{2s-1}
\left[\frac{(-)^t}{n^t}+\frac{1}{(n+1)^t}\right].
\end{equation}

Sums of reciprocal powers of the even or odd integers are
in terms of Riemann's $\zeta$-function
\cite[0.233]{GR}\cite[(335)]{Jolley}
\begin{equation}
\sum_{k=1}^\infty \frac{1}{(2k)^t}
=
\frac{1}{2^t}\zeta(t),
\end{equation}
and
\begin{equation}
\sum_{k=1}^\infty \frac{1}{(2k+1)^t}
= 
[1-\frac{1}{2^t}]\zeta(t) -1.
\end{equation}
Combining the previous three equations generates
(with a little extra care at $t=1$ \cite[0.234]{GR})
\begin{eqnarray}
&&T_e(2s)\equiv \sum_{k=1}^\infty \frac{1}{(2k)^{2s}(2k+1)^{2s}}
\nonumber\\
&&=
\sum_{k=1}^\infty \sum_{t=1}^{2s}
\binom{4s-t-1}{2s-1}
\left[
\frac{(-)^t}{(2k)^t}+\frac{1}{(2k+1)^t}
\right]
\nonumber\\
&&=
\sum_{k=1}^\infty \bigg\{
-\frac{\binom{4s-2}{2s-1}}{(2k)(2k+1)}
\nonumber\\
&&
\qquad +\sum_{t=2}^{2s}
\binom{4s-t-1}{2s-1}
\left[
\frac{(-)^t}{(2k)^t}+\frac{1}{(2k+1)^t}
\right]
\bigg\}
\nonumber\\
&&=
-\binom{4s-2}{2s-1}[1-\ln 2]
\nonumber\\
&&
\qquad
+\sum_{k=1}^\infty \sum_{t=2}^{2s}
\binom{4s-t-1}{2s-1}
\left[
\frac{(-)^t}{(2k)^t}+\frac{1}{(2k+1)^t}
\right]
\nonumber\\
&&=
\sum_{t=1}^{2s}
\binom{4s-t-1}{2s-1}
\left\{
-1 +[1-\frac{1-(-)^t}{2^t}]\zeta(t)
\right\}
.
\label{eq.ksum}
\end{eqnarray}
For odd $t$, $[1-\frac{1-(-)^t}{2^t}]\zeta(t)$ equals Dirichlet's $\eta$-function,
in particular $\eta(1)=\log 2$ at the pole of $\zeta(1)$ \cite[Tab.\ 23.3]{AS}.

The three base examples of this format are:
\begin{eqnarray}
&&\sum_{k=1}^\infty
\frac{1}{(2k)^2(2k+1)^2}
=-3+\frac{\pi^2}{6}+2\log 2
\nonumber\\
&&\approx 0.03122842796811705530687941 ;
\end{eqnarray}
\begin{eqnarray}
&&\sum_{k=1}^\infty
\frac{1}{(2k)^4(2k+1)^4}
=-35+\frac{\pi^4}{90}+3\zeta(3)+\frac{5\pi^2}{3}+20\log 2
\nonumber\\
&&\approx 0.00077822287109160078401223 ;
\end{eqnarray}
\begin{eqnarray}
&&\sum_{k=1}^\infty
\frac{1}{(2k)^6(2k+1)^6}
=-462+21\pi^2 +42\zeta(3)
\nonumber\\
&&
+252 \log 2
+\frac{7\pi^4}{30}+\frac{\pi^6}{945}
+\frac{45}{8}\zeta(5)
\nonumber\\
&&\approx 0.0000214492855159526203348.
\end{eqnarray}
Interchange of the summation over the Taylor orders and over the $k$ and
insertion of (\ref{eq.ksum}) into (\ref{eq.logcost}) leads to the value (\ref{eq.numlog}).
Exponentiation gives
\begin{eqnarray}
C_e \equiv
\prod_{n=2,4,6,\ldots}^\infty \cos\frac{\pi}{n(n+1)}
\nonumber\\
\approx 0.85135730526671405636170.
\label{eq.Ce}
\end{eqnarray}

\subsection{Odd lower term in the product}
In (\ref{eq.numlog}), the smaller factor in the product $n(n+1)$ was always even.
With exactly the same technique we obtain a ``complete'' version of
(\ref{eq.ksum}) where the smaller term steps through all positive integers:
\begin{eqnarray}
&&T(2s)\equiv \sum_{n=1}^\infty \frac{1}{n^{2s}(n+1)^{2s}}
\nonumber\\
&&=
\sum_{t=1}^{2s}\binom{4s-t-1}{2s-1}
\left\{
[1+(-)^t]\zeta(t)-1
\right\}.
\end{eqnarray}
Here $[1+(-1)^t]\zeta(t)$ is to be interpreted as $0$ if $t=1$,
ignoring the pole of $\zeta$.
The three basic examples are
\begin{eqnarray}
&&\sum_{n=1}^\infty \frac{1}{n^2(n+1)^2} = \frac{\pi^2}{3}-3
\nonumber\\
&&\approx 0.2898681336964528729448303333 ;
\end{eqnarray}
\begin{eqnarray}
&&\sum_{n=1}^\infty \frac{1}{n^4(n+1)^4} = -35+\frac{10\pi^2}{3}+\frac{\pi^4}{45}
\nonumber\\
&&\approx 0.0633278043868051124803107260 ;
\end{eqnarray}
\begin{eqnarray}
&&\sum_{n=1}^\infty \frac{1}{n^6(n+1)^6} = -462+42\pi^2+\frac{7\pi^4}{15}+\frac{2\pi^6}{945}
\nonumber\\
&&\approx 0.0156467855897643141498131091.
\end{eqnarray}
The complete version of (\ref{eq.numlog}) does not exist because
the term at $n=1$ contributes
$\log \cos (\pi/2) = -\infty$.
We drop this term at $n=1$ and use
$T(2s)-2^{-2s}$ with (\ref{eq.taylogcos}) to compute
\begin{equation}
\log \prod_{n=2}^\infty \cos \frac{\pi}{n(n+1)}
\approx
-0.2039684236116246918364049,
\end{equation}
and its exponential value
\begin{equation}
C\equiv \prod_{n=2}^\infty \cos \frac{\pi}{n(n+1)}
\approx
0.815488120950370848344387.
\label{eq.C}
\end{equation}

If the smaller factor in $n(n+1)$ is odd, the difference is involved:
\begin{equation}
T_o(2s) \equiv \sum_{n=1,3,5,\ldots}^\infty \frac{1}{n^{2s}(n+1)^{2s}}
=T(2s)-T_e(2s).
\end{equation}
Division of (\ref{eq.C}) through (\ref{eq.Ce}) yields the complement
\begin{eqnarray}
&& C_o \equiv \prod_{n=3,5,7\ldots}^\infty \cos \frac{\pi}{n(n+1)}
=\frac{C}{C_e}
\nonumber \\
&& \approx
0.95786823687957188013580826171688
\label{eq.Co}
\end{eqnarray}
for use in (\ref{eq.Couse}).

\subsection{Numerator $2\pi$ with odd lower term}
The factor $2\pi$ in the numerator of (\ref{eq.roteoInf}),
\begin{equation}
\prod_{k=3,5,7\ldots}\cos\frac{2\pi}{k(k+1)}
=
\prod_{n=1}^\infty \cos\frac{\pi}{(2n+1)(n+1)}
\end{equation}
causes slow convergence of the methods shown above.
An acceleration method with deferred summation is available \cite{SebahGourdon}:
The partial product up to some $n\le M$ is calculated
explicitly, and the logarithm
of the remaining infinite product is expanded in a Taylor series in $1/n$:
\begin{eqnarray}
&& \sum_{n>M} \log \cos\frac{\pi}{(2n+1)(n+1)}
\nonumber\\
&&
=\sum_{n>M} [
-\frac{\pi^2}{8}\frac{1}{n^4}
+\frac{3\pi^2}{8}\frac{1}{n^5}
-\frac{23\pi^2}{32}\frac{1}{n^6}
+\frac{9\pi^2}{8}\frac{1}{n^7}+\cdots
]
\label{eq.cos2pT}
\end{eqnarray}
Each term on the right hand side is then replaced by an incomplete $\zeta$-function,
\begin{equation}
\sum_{n>M} \frac{1}{n^s} = \zeta(s)-\sum_{n=1}^M\frac{1}{n^s}.
\end{equation}
With $M=10$ and (\ref{eq.cos2pT}) expanded
up to $O(n^{-30})$
we obtain for example
\begin{equation}
\prod_{k=3,5,7\ldots}\cos\frac{2\pi}{k(k+1)}
\approx
0.8373758680415481080004775.
\label{eq.cos2pi}
\end{equation}

\section{Quenching Factor of Kitson's  Constant}\label{sec.logp}

The logarithm of the product in
(\ref{eq.ratPri})
is the a sum over all odd primes $p_j\ge 3$:
\begin{equation}
\log \prod_{p_j\ge 3} \cos\frac{\pi}{p_jp_{j+1}}
=
\sum_{p_j\ge 3} \log \cos\frac{\pi}{p_jp_{j+1}}
\end{equation}
Again with (\ref{eq.taylogcos})
we evaluate
\begin{equation}
T(2s)\equiv \sum_{p_j} \frac{1}{p_j^{2s}p_{j+1}^{2s}}
\end{equation}
for integer $s$.
\begin{eqnarray}
T(2)&\approx & 0.005519522774559;
\\
T(4)&\approx & 0.0000204508599535;
\\
T(6)&\approx & 0.000000088340410739027;
\\
T(8)&\approx & 0.000000000390629312549651477,
\end{eqnarray}
so
\begin{equation}
\sum_{p_j\ge 3} \log\cos\frac{\pi}{p_jp_{j+1}} \approx -0.02740567
\end{equation}
and after exponentiation
\begin{eqnarray}
\prod_{p_j\ge 3} \cos\frac{\pi}{p_jp_{j+1}}
\approx 0.9729664541346255360938192\ldots
\label{eq.Cprim}
\end{eqnarray}

\bibliography{all}

\begin{thebibliography}{17}%
\makeatletter
\providecommand \@ifxundefined [1]{%
 \@ifx{#1\undefined}
}%
\providecommand \@ifnum [1]{%
 \ifnum #1\expandafter \@firstoftwo
 \else \expandafter \@secondoftwo
 \fi
}%
\providecommand \@ifx [1]{%
 \ifx #1\expandafter \@firstoftwo
 \else \expandafter \@secondoftwo
 \fi
}%
\providecommand \natexlab [1]{#1}%
\providecommand \enquote  [1]{``#1''}%
\providecommand \bibnamefont  [1]{#1}%
\providecommand \bibfnamefont [1]{#1}%
\providecommand \citenamefont [1]{#1}%
\providecommand \href@noop [0]{\@secondoftwo}%
\providecommand \href [0]{\begingroup \@sanitize@url \@href}%
\providecommand \@href[1]{\@@startlink{#1}\@@href}%
\providecommand \@@href[1]{\endgroup#1\@@endlink}%
\providecommand \@sanitize@url [0]{\catcode `\\12\catcode `\$12\catcode
  `\&12\catcode `\#12\catcode `\^12\catcode `\_12\catcode `\%12\relax}%
\providecommand \@@startlink[1]{}%
\providecommand \@@endlink[0]{}%
\providecommand \url  [0]{\begingroup\@sanitize@url \@url }%
\providecommand \@url [1]{\endgroup\@href {#1}{\urlprefix }}%
\providecommand \urlprefix  [0]{URL }%
\providecommand \Eprint [0]{\href }%
\providecommand \doibase [0]{http://dx.doi.org/}%
\providecommand \selectlanguage [0]{\@gobble}%
\providecommand \bibinfo  [0]{\@secondoftwo}%
\providecommand \bibfield  [0]{\@secondoftwo}%
\providecommand \translation [1]{[#1]}%
\providecommand \BibitemOpen [0]{}%
\providecommand \bibitemStop [0]{}%
\providecommand \bibitemNoStop [0]{.\EOS\space}%
\providecommand \EOS [0]{\spacefactor3000\relax}%
\providecommand \BibitemShut  [1]{\csname bibitem#1\endcsname}%
\let\auto@bib@innerbib\@empty
\bibitem [{\citenamefont {Finch}(2003)}]{Finch}%
  \BibitemOpen
  \bibfield  {author} {\bibinfo {author} {\bibfnamefont {S.~R.}\ \bibnamefont
  {Finch}},\ }\href@noop {} {\emph {\bibinfo {title} {Mathematical
  Constants}}},\ \bibinfo {series} {Encyclopedia of Mathematics and its
  Applications}\ No.~\bibinfo {number} {94}\ (\bibinfo  {publisher} {Cambridge
  University Press},\ \bibinfo {address} {Cambridge},\ \bibinfo {year}
  {2003})\BibitemShut {NoStop}%
\bibitem [{\citenamefont {Weisstein}(2002)}]{WeissteinCRC}%
  \BibitemOpen
  \bibinfo {editor} {\bibfnamefont {E.~W.}\ \bibnamefont {Weisstein}},\ ed.,\
  \href@noop {} {\emph {\bibinfo {title} {{CRC} Concise Encyclopedia of
  Mathematics}}},\ \bibinfo {edition} {2nd}\ ed.\ (\bibinfo  {publisher}
  {Chapman \& Hall/CRC},\ \bibinfo {year} {2002})\BibitemShut {NoStop}%
\bibitem [{\citenamefont {Sloane}(2003)}]{EIS}%
  \BibitemOpen
  \bibfield  {author} {\bibinfo {author} {\bibfnamefont {N.~J.~A.}\
  \bibnamefont {Sloane}},\ }\href {http://oeis.org/} {\bibfield  {journal}
  {\bibinfo  {journal} {Notices Am.\ Math.\ Soc.}\ }\textbf {\bibinfo {volume}
  {50}},\ \bibinfo {pages} {912} (\bibinfo {year} {2003})},\ \bibinfo {note}
  {http://oeis.org/},\ \Eprint {http://arxiv.org/abs/math.CO/0312448}
  {arXiv:math.CO/0312448} \BibitemShut {NoStop}%
\bibitem [{\citenamefont {Abramowitz}\ and\ \citenamefont {Stegun}(1972)}]{AS}%
  \BibitemOpen
  \bibinfo {editor} {\bibfnamefont {M.}~\bibnamefont {Abramowitz}}\ and\
  \bibinfo {editor} {\bibfnamefont {I.~A.}\ \bibnamefont {Stegun}},\ eds.,\
  \href@noop {} {\emph {\bibinfo {title} {Handbook of Mathematical
  Functions}}},\ \bibinfo {edition} {9th}\ ed.\ (\bibinfo  {publisher} {Dover
  Publications},\ \bibinfo {address} {New York},\ \bibinfo {year}
  {1972})\BibitemShut {NoStop}%
\bibitem [{\citenamefont {Conway}\ \emph {et~al.}(1998)\citenamefont {Conway},
  \citenamefont {Radin},\ and\ \citenamefont {Sadun}}]{Conwayarxiv98}%
  \BibitemOpen
  \bibfield  {author} {\bibinfo {author} {\bibfnamefont {J.~H.}\ \bibnamefont
  {Conway}}, \bibinfo {author} {\bibfnamefont {C.}~\bibnamefont {Radin}}, \
  and\ \bibinfo {author} {\bibfnamefont {L.}~\bibnamefont {Sadun}},\
  }\href@noop {} {\bibfield  {journal} {\bibinfo  {journal}
  {arXiv:math-ph/9812019}\ } (\bibinfo {year} {1998})}\BibitemShut {NoStop}%
\bibitem [{\citenamefont {Girstmair}(1997)}]{GirstmairAA81}%
  \BibitemOpen
  \bibfield  {author} {\bibinfo {author} {\bibfnamefont {K.}~\bibnamefont
  {Girstmair}},\ }\href@noop {} {\bibfield  {journal} {\bibinfo  {journal}
  {Acta Arithm.}\ }\textbf {\bibinfo {volume} {81}},\ \bibinfo {pages} {387}
  (\bibinfo {year} {1997})}\BibitemShut {NoStop}%
\bibitem [{\citenamefont {Grimstone}(1980)}]{GrimstoneMG64}%
  \BibitemOpen
  \bibfield  {author} {\bibinfo {author} {\bibfnamefont {C.~J.}\ \bibnamefont
  {Grimstone}},\ }\href@noop {} {\bibfield  {journal} {\bibinfo  {journal}
  {Math. Gaz.}\ }\textbf {\bibinfo {volume} {64}},\ \bibinfo {pages} {120}
  (\bibinfo {year} {1980})}\BibitemShut {NoStop}%
\bibitem [{\citenamefont {Kitson}(2006)}]{KitsonArxiv06}%
  \BibitemOpen
  \bibfield  {author} {\bibinfo {author} {\bibfnamefont {A.~R.}\ \bibnamefont
  {Kitson}},\ }\href@noop {} {\bibfield  {journal} {\bibinfo  {journal}
  {arXiv:math.HO/0608186}\ } (\bibinfo {year} {2006})}\BibitemShut {NoStop}%
\bibitem [{\citenamefont {Kitson}(2008)}]{KitsonMG92}%
  \BibitemOpen
  \bibfield  {author} {\bibinfo {author} {\bibfnamefont {A.~R.}\ \bibnamefont
  {Kitson}},\ }\href@noop {} {\bibfield  {journal} {\bibinfo  {journal} {Math.
  Gaz}\ }\textbf {\bibinfo {volume} {92}},\ \bibinfo {pages} {293} (\bibinfo
  {year} {2008})}\BibitemShut {NoStop}%
\bibitem [{\citenamefont {Stephens}(1995)}]{StephensMG79}%
  \BibitemOpen
  \bibfield  {author} {\bibinfo {author} {\bibfnamefont {E.}~\bibnamefont
  {Stephens}},\ }\href@noop {} {\bibfield  {journal} {\bibinfo  {journal}
  {Math. Gaz.}\ }\textbf {\bibinfo {volume} {79}},\ \bibinfo {pages} {561}
  (\bibinfo {year} {1995})}\BibitemShut {NoStop}%
\bibitem [{\citenamefont {Gradstein}\ and\ \citenamefont {Ryshik}(1981)}]{GR}%
  \BibitemOpen
  \bibfield  {author} {\bibinfo {author} {\bibfnamefont {I.}~\bibnamefont
  {Gradstein}}\ and\ \bibinfo {author} {\bibfnamefont {I.}~\bibnamefont
  {Ryshik}},\ }\href@noop {} {\emph {\bibinfo {title} {Summen-, {P}rodukt- und
  {I}ntegraltafeln}}},\ \bibinfo {edition} {1st}\ ed.\ (\bibinfo  {publisher}
  {Harri Deutsch},\ \bibinfo {address} {Thun},\ \bibinfo {year}
  {1981})\BibitemShut {NoStop}%
\bibitem [{\citenamefont {Mahoney}\ and\ \citenamefont
  {Sivazlian}(1983)}]{MahoneyJCAM9}%
  \BibitemOpen
  \bibfield  {author} {\bibinfo {author} {\bibfnamefont {J.~J.}\ \bibnamefont
  {Mahoney}}\ and\ \bibinfo {author} {\bibfnamefont {B.~D.}\ \bibnamefont
  {Sivazlian}},\ }\href {\doibase 10.1016/0377-0427(83)90018-3} {\bibfield
  {journal} {\bibinfo  {journal} {J. Comp. Appl. Math.}\ }\textbf {\bibinfo
  {volume} {9}},\ \bibinfo {pages} {247} (\bibinfo {year} {1983})}\BibitemShut
  {NoStop}%
\bibitem [{\citenamefont {Velleman}(2002)}]{VellemanAMM109}%
  \BibitemOpen
  \bibfield  {author} {\bibinfo {author} {\bibfnamefont {D.~J.}\ \bibnamefont
  {Velleman}},\ }\href@noop {} {\bibfield  {journal} {\bibinfo  {journal} {Am.
  Math. Monthly}\ }\textbf {\bibinfo {volume} {109}},\ \bibinfo {pages} {746}
  (\bibinfo {year} {2002})}\BibitemShut {NoStop}%
\bibitem [{\citenamefont {Xin}(2004)}]{Xinarxiv04}%
  \BibitemOpen
  \bibfield  {author} {\bibinfo {author} {\bibfnamefont {G.}~\bibnamefont
  {Xin}},\ }\href {math/0408189} {\bibfield  {journal} {\bibinfo  {journal}
  {arXiv}\ } (\bibinfo {year} {2004})},\ \Eprint
  {http://arxiv.org/abs/math/0408189} {arXiv:math/0408189} \BibitemShut
  {NoStop}%
\bibitem [{\citenamefont {Eustice}\ and\ \citenamefont
  {Klamkin}(1979)}]{EusticeAMM86}%
  \BibitemOpen
  \bibfield  {author} {\bibinfo {author} {\bibfnamefont {D.}~\bibnamefont
  {Eustice}}\ and\ \bibinfo {author} {\bibfnamefont {M.~S.}\ \bibnamefont
  {Klamkin}},\ }\href@noop {} {\bibfield  {journal} {\bibinfo  {journal} {Am.
  Math. Monthly}\ }\textbf {\bibinfo {volume} {86}},\ \bibinfo {pages} {478}
  (\bibinfo {year} {1979})}\BibitemShut {NoStop}%
\bibitem [{\citenamefont {Jolley}(1961)}]{Jolley}%
  \BibitemOpen
  \bibinfo {editor} {\bibfnamefont {J.~B.~W.}\ \bibnamefont {Jolley}},\ ed.,\
  \href@noop {} {\emph {\bibinfo {title} {Summation of series}}},\ \bibinfo
  {edition} {2nd}\ ed.\ (\bibinfo  {publisher} {Dover Publications},\ \bibinfo
  {address} {New York},\ \bibinfo {year} {1961})\BibitemShut {NoStop}%
\bibitem [{\citenamefont {Sebah}\ and\ \citenamefont
  {Gourdon}(2001)}]{SebahGourdon}%
  \BibitemOpen
  \bibfield  {author} {\bibinfo {author} {\bibfnamefont {P.}~\bibnamefont
  {Sebah}}\ and\ \bibinfo {author} {\bibfnamefont {X.}~\bibnamefont
  {Gourdon}},\ }\href
  {http://numbers.computation.free.fr/Constants/constants.html} {\enquote
  {\bibinfo {title} {Constants from number theory},}\ } (\bibinfo {year}
  {2001}),\ \bibinfo {note}
  {http://numbers.computation.free.fr/Constants/constants.html}\BibitemShut
  {NoStop}%
\end{thebibliography}%

\end{document}